\documentclass[3p]{elsarticle}
\usepackage{epsfig}
\usepackage{amsfonts}
\usepackage{amssymb}
\usepackage{graphics}
\usepackage{mathrsfs}
\usepackage{amsmath}
\usepackage{mathtools}
\usepackage{mathrsfs}
\usepackage{float}
\usepackage{setspace}
\newtheorem{theorem}{Theorem}[section]

\newtheorem{example}{Example}[section]

\newtheorem{definition}{Definition}[section]

\newproof{pol}{\bf{Proof of Lemma}}
\newproof{pot}{\bf{Proof}}
\newproof{pf}{\bf{Proof of Theorem}}

\usepackage{epstopdf}
\usepackage{algorithm}
\usepackage{algpseudocode}
\usepackage{tabularx,ragged2e,booktabs,caption}

\journal{Applied Mathematics and Computation}
\begin{document}
	
	
	\begin{frontmatter}
		
		\title{\bf Unpredictable Strings}
		
		\author{Marat Akhmet$^{*,a}$} \ead{marat@metu.edu.tr}
		\author{Astrit Tola$^{a}$} \ead{astrit.tola@gmail.com}
		
		\address{$^{a}$Department of Mathematics, Middle East Technical University, 06800 Ankara, Turkey}
		\address{$^*$Corresponding Author Tel.:+90 312 210 5355}
		
		\begin{abstract} 
			
	A novel notion of unpredictable strings is revealed and utilized to define deterministic unpredictable sequences on a finite number of symbols. We prove the first law of large strings for random processes in discrete time, which confirms that there exists the uncountable set of unpredictable realizations. The hypothesis on the second law of large strings is formulated, which is relative to the Bernoulli theorem. Theoretical and numerical backgrounds for the   phenomenon are provided.
			
		\end{abstract}
		
		
	\end{frontmatter}
	
	\section{Introduction}
  
  We  have    developed   the concept  of   strong  relation between   deterministic  chaos  and random dynamics  in our  recent papers  \cite{Akh51}-\cite{Akh53}.  This time,  the   notions of  unpredictable strings of symbols and infinite unpredictable sequences with unpredictable strings of unbounded lengths are  introduced. The definitions strongly relate to the concept of the unpredictable point  \cite{Akh16a}-\cite{Akh19S}. A numerical simulation of the Bernoulli process is performed to demonstrate that the realization is a part of an unpredictable sequence. A numerical analysis
  confirms that specific properties for the random dynamics are valid, the first and second laws of large (unpredictable)  strings. Besides, a Matlab algorithm to construct sequences with inductively increasing lengths of unpredictable strings is provided.

	\section{The unpredictable strings}\label{Intro}
\begin{figure}[ht]
	\centering
	\includegraphics[height=4cm]{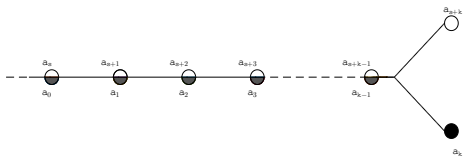}
	\caption{The illustration of the unpredictable string of length $k$.}
	\label{fig2.3}
	\end{figure}

	In this section, we introduce the main concept of this paper, unpredictable strings, and utilize them to determine unpredictable sequences. Let $a_i$, $i=0,1,2,...$, be an infinite sequence of symbols. The diagram in Figure \ref{fig2.3} illustrates the definition.
	\begin{definition}
		A finite array $(a_s,a_{s+1},...,a_{s+k})$, where $s$ and $k$ are positive integers, is said to be an unpredictable string of length $k$ if $a_i=a_{s+i}$, for $i=0,1,2,...,k-1$, and $a_k \neq a_{s+k}$.
	\end{definition}

	\begin{definition}
	\label{01}
	The sequence $a_i$ is unpredictable if  it  admits unpredictable strings with arbitrary large lengths.
	\end{definition}	
\begin{definition}\label{02} 
		\cite{Akh18}
		The sequence $a_i$ is unpredictable if there exist sequences $\zeta _n$, $\eta_n$ of positive integers both of which diverge to infinity such that $a_{\zeta_n+l} = a_l$, $l=0,1,2,...,\eta_n-1$, and $a_{\zeta_n+\eta_n} \neq a_{\eta_n}$, for each $n \in \mathbb{N}$.
	\end{definition}
	    \begin{theorem}
	\label{03}
	The  Definitions \ref{01} and \ref{02} are  equivalent.
    \end{theorem}
	\begin{pot}
		Let sequence $a_i$ be unpredictable. Then the finite arrays $(a_{\zeta_n},a_{\zeta_n+1},...,a_{\zeta_n+\eta_n})$ are unpredictable strings of length $\eta_n$, for each natural $n$. Thus the sequence admits unpredictable strings with arbitrary large lengths.

		Conversely let $a_i$ be a sequence that admits unpredictable strings of arbitrary large lengths, i. e., there is a sequence $i_n$, $n=1,2,3,...$, such that the finite arrays $(a_{i_n},a_{i_n+1},...,a_{i_n+k})$ are unpredictable strings. By setting $\zeta _n=i_n$ and $\eta_n=i_n+k$ we deduce that the sequence $a_i$ is unpredictable in  light of Definition \ref{02}. $\Box$
	\end{pot}

Fix a positive integer $k$ and denote by $S_k$  the sets  of all indexes $s$ such that the strings $(a_s,a_{s+1},...,a_{s+k})$ are unpredictable within   the sequence $a_i, i = 1, 2, ... $, which is not necessarily unpredictable.

\begin{theorem}
	\label{1}
	The sets $S_l$ and $S_q$ do not intersect if $l < q$.
\end{theorem}
\begin{pot}
	Assume,  on contrary, that   sets $S_l$ and $S_q$  have a common element  $s.$  Then, we have that $a_l\neq a_{s+l}$ if $s\in S_l$ and $a_{l}=a_{s+l}$ if $s\in S_q$.   This  contradiction    completes the prove. $\Box$	
\end{pot}

\begin{theorem}
	\label{2}
	Assume that $a_i$ is an unpredictable sequence. Then each $a_j$ with positive $j$ is the first element of an unpredictable string, if $a_j = a_0$.
\end{theorem}
\begin{pot}
	Assume the opposite. Then  one can  show that  the sequence $\alpha$  is periodic one.  That  is not  unpredictable sequence.  $\Box$	
\end{pot}

	\section{Numerical analysis of the Bernoulli process} \label{Bern}

In this section, we will scrutinize the realizations of Bernoulli processes, by considering them as sequences consisting of the digits 1 and 0 with positive probabilities.

First, we  will  build unpredictable strings of inductively increasing lengths by using fixed complex vectors, $v_1,v_2,...,v_r$. 

Let us set  $a_0=random(\{v_1,v_2,...,v_r\})$ and $a_1=random(\{v_1,v_2,...,v_r\})$. Then for increasing $k=1,2,3,...$, we define $a_{m(k)+j}=a_j$, for $j<k$, and $a_{m(k)+j}=random(\{v_1,v_2,...,v_r\}-a_j)$, for $j=k$, where $m(k+1)=m(k)+k$ with $m(1)=2$. 

 The immediately following Algorithm \ref{Algseq} is the scalar case for $r=2$, $v_1=0$ and $v_2=1$.
The sequence $(0, 1, 0, 0, 0, 1, 1, 0, 1, 0, 1, 0, 1, 0, 0, 1, 0, 1, 0, 0, 0, 0, 0, 1, 0, 0, 0, 1, 0, ...)$ is a result of the algorithm application

\begin{algorithm}[h!]
	\caption{\,\,\,\,Unpredictable sequences}
	\label{Algseq}
	\begin{algorithmic}[1]
		
		\State {$m=2$}
		\For {$k=1,2,3,...$}
		\State {$a_0=0$}
		\State {$a_1=1$}
		\For {$j=0:k$}
		\If {$j<k$} 
		\State {$a_{m+j}=a_j$}
		\ElsIf {$j=k$}
		\State {$a_{m+j}\neq a_j$}
		\State {$m=m+k$}
		\EndIf
		\EndFor
		\EndFor
	\end{algorithmic}
\end{algorithm}

Let us introduce several characteristics that are of usage for analyses of finite realizations. For fixed natural number $m$, consider a finite realization, $a_i, i = 0, 1, ..., m$. Denote by $K(m)$ the largest length of unpredictable strings in the array. For every $k$ between $1$ and $K(m)$, denote by $q_k$  the number of $k-$lengthy unpredictable strings within the array, by $\xi_k$ the largest index such that $(a_{\xi_k}, a_{\xi_{k+1}}, ..., a_{\xi_k+k})$ is an
unpredictable string within the array, and by $N(m)$ the number of all unpredictable strings, which have a non-empty intersection with the array.

Next, we provide statistical results on the realization, which are obtained by Matlab simulations for the Bernoulli process with probability $p = 0.6$ and $m = 9 \times 10^5$. We have evaluated values of $K(n)$, $\xi_{K(n)}$ and $N(n)/n$, for each $n$ from $1$ to $m$. Ten samples of the simulations are provided in Table \ref{tab7}. According to the
full data obtained in simulations, the realization can be considered as part of an unpredictable sequence, since there are unpredictable strings with increasing lengths. Moreover, $N(n)/n \approx p$, if $n$ is large.

	\begin{table}[ht]
		\centering
		\begin{tabular}{c  c  c  c}
			\hline
			$n$ &$K(n)$&$\xi_{K(n)}$ &$ N(n)/n$\\
			\hline
			50		&10	& 20 & 0.72\\
			
			
			200		&10 & 20&0.58\\
			
			500		&10 & 228&0.586\\
			
			
			2000	&14	&1008&0.596\\
			
			
			10000	&14	&3469&0.6031\\
			
			20000	&18	&19206&0.5995\\
			
			
			
			100000	&21	&74683&0.6014\\
			
			
			
			
			500000	&21	&401088&0.6003\\
			
			
			
			800000	&21	&663684&0.6001\\
			
			900000	&28	&874766&0.5686\\
			\hline
			
		\end{tabular}
		\caption{The values ${K(n)}$, $\xi_{K(n)}$ and $N(n)/n$ for the finite realization.}
		\label{tab7}
	\end{table}

	\section{Laws of large strings}
In this section, we consider a discrete-time random process $\textbf{X}(n)$ with  the finite state space of $r$ different symbols  $s_1, s_2,..., s_r.$ 
The function admits values $s_i$ with positive probabilities $p_i$, $i=1,2,...,r$, which sum is equal to the unit. A realization $\alpha$ of the process is the  sequence $a_i$, $i=1,2,...$, and  a finite realization $\alpha_m$ is the  array $a_i$, $i=1,2,...,m$. We claim that stochastic processes with discrete-timeand finite-state spaces satisfy the following theorem.
\begin{theorem}
(the first law of large strings) The discrete time random process $\textbf{X}(n)$ with the finite state space admits uncountable set of realizations, which are unpredictable sequences in the sense of Definition \ref{01}.
\end{theorem}
\begin{pot}
Let us consider the space $\Sigma_r $   of infinite sequences  of  finite set of symbols $s_1, s_2,..., s_r,$   with  the metric
\begin{eqnarray}
d(\xi,\zeta)={\Sigma}_{k=0}^\infty\dfrac{|\xi_k-\zeta_k|}{2^k},
\end{eqnarray}
where $\xi = (\xi_0\xi_1\xi_2...), \zeta = (\zeta_0\zeta_1\zeta_2...)$. The Bernoulli shift $\sigma$    on  $\Sigma_r$    is defined as $\sigma(\xi_0\xi_1\xi_2...) =(\xi_1\xi_2\xi_3...)$. The map  is continuous and $\Sigma_r$ is a compact metric space \cite{Wiggins}.

It is clear that the set of all realizations of the random dynamics $X\textbf{}(n)$ coincides with the set of all sequences of the symbolic dynamics on  $\Sigma_r.$   According to the result in \cite{Akh16b}, the symbolic dynamics admits an unpredictable point, $i^*$, a sequence from the set $\Sigma_r.$    There  is the  uncountable  set  of   unpredictable  points, which  are unpredictable sequences  in the  sense of Definition \ref{01}. 
$\Box$
\end{pot}
  It is important that the  set   is  the  closure for the unpredictable orbit. The density   is  considered in the shift dynamics sense.  The property   of the metric implies that each arc of any sequence in the space coincides with some arc of the unpredictable sequence.  
	
Denote by n(m) the number of elements, which are equal to $a_0$ in a finite string. The limit $E[a_0] =\lim_{m\to\infty} n(m)/m$ is said to be the expected value \cite{Casta}. It is clear that $E[a_0] = p_i$, if $a_0 = s_i, i= 1, ..., r$.
 
 Theorem \ref{2} implies the equality N(m) = n(m). Hence, by the Bernoulli theorem and arguments
 assumed for the first law, one may suggest that the following second law of large strings may be valid.

\begin{theorem}
	If the discrete time random process $X(n)$ admits a finite state space, then the relation
	\begin{eqnarray}
	\label{p2}
	\lim_{m\to\infty}P \left( \left|\dfrac{N(m)}{m}-E[a_0] \right|<\varepsilon \right) =1
    \end{eqnarray}
	holds for any $\varepsilon> 0$.
	
\end{theorem}

We can not prove the last theorem yet rigorously. This is why, we suggest it as an open problem.
Since of Theorem 2.3, the following assertion is correct, which can be useful for applications.	
	
	\begin{theorem}
		\label{4}  If a realization $\alpha$ is an unpredictable sequence, then the relation
		\begin{eqnarray}
		\label{p1}
		\lim_{m\to\infty}\dfrac{N(m)}{m} =E[a_0]
		\end{eqnarray}
		is valid.
	\end{theorem}

\begin{example}  To    have  more impressions  of  the unpredictable  strings,   let  us  consider  the    graph  of the    piece-wise constant  function,  $H(t),$  which    values on intervals      $[i/10,(i+1)/10),    i =0,1,\ldots,199,$  are  assigned   randomly  $1$  or  $-1$    with  equal  probability $1/2.$    The   two  unpredictable  strings   as result  of the  Bernoulli  process  are present,  in the  red, in the  Figure \ref{fig3.3},  $(a).$   The second one,   with  length    of  $0.7$  units, placed  between   coordinates   $14$   and $16,$  shown  in    Figure \ref{fig3.3},  $(c),$    while  its  corresponding   initial   arc,  in  Figure \ref{fig3.3},  $(b).$
	
	\begin{figure}[ht]
		\centering
		\includegraphics[height=6cm]{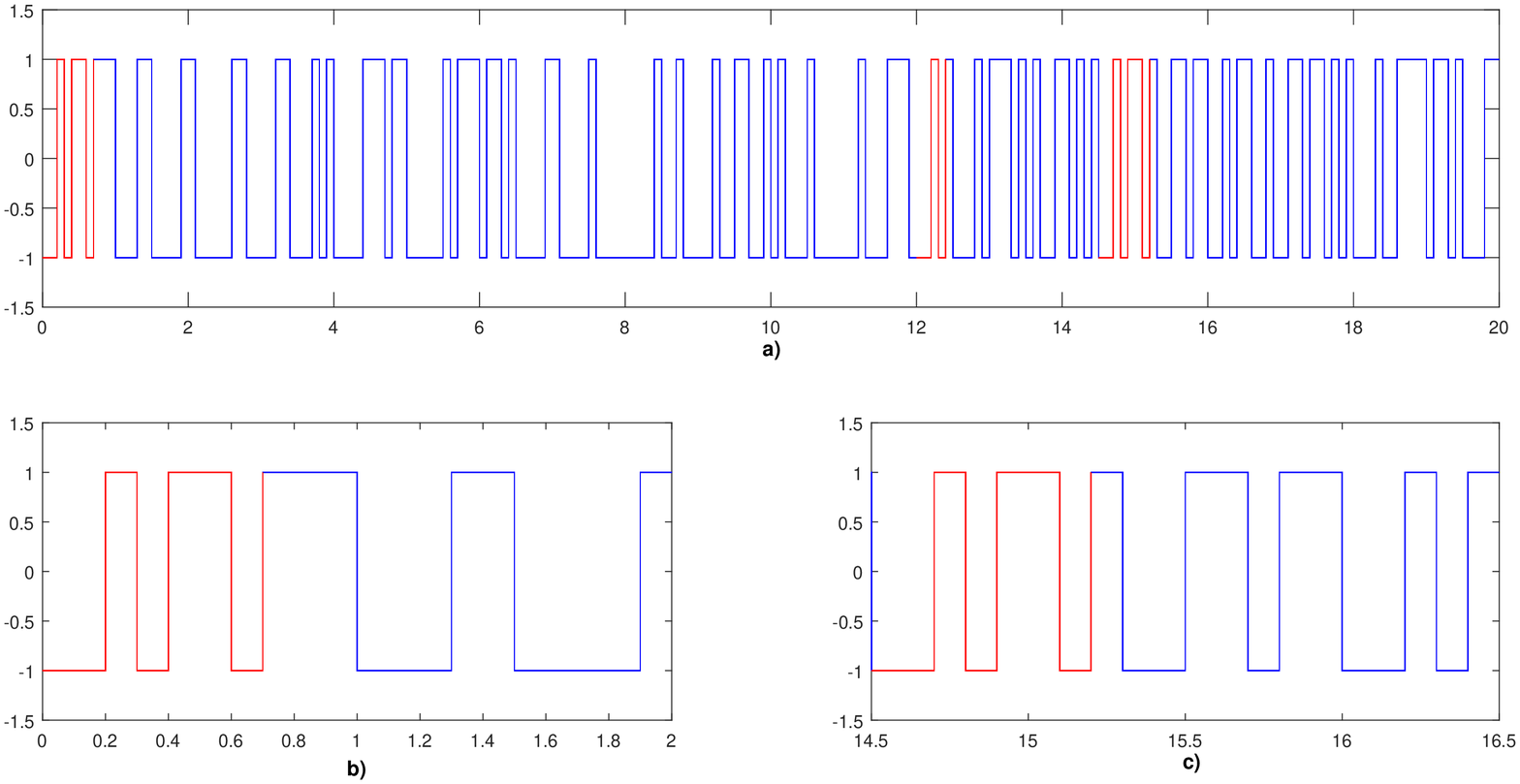}
		\caption{The graph  of the  function $H(t),$  which  illustrates   unpredictable  strings  appearance.}
		\label{fig3.3}
	\end{figure}
	
\end{example}	

%
%
%


\section*{References}
	
	\begin{thebibliography}{99}
		
		\bibitem{Akh51}  {\textsc{Akhmet, M.},  Modular chaos for random processes, 2020,  arXiv:2004.08383.}		
		
			\bibitem{Akh52}  {\textsc{Akhmet, M.},  Domain structured chaos for discrete random processes, 2020, arXiv:1912.10478 }
				
			\bibitem{Akh53} 	{\textsc{Akhmet, M., Fen, M. O., Alejaily, E. M.},  A randomly determined unpredictable function,   2019, arXiv:1910.12758 .}
			
			\bibitem{Akh16a} {\textsc{Akhmet, M., Fen, M. O.}, Unpredictable points and chaos, \textit{Commun. Nonlinear Sci. Numer. Simulat.} \textbf{40} (2016), 1--5.}
			
			\bibitem{Akh16b} {\textsc{Akhmet, M., Fen, M. O.}, Poincar\'{e} chaos and unpredictable functions, \textit{Commun. Nonlinear Sci. Numer. Simulat.} \textbf{48} (2016), 85--94.}
		
		\bibitem{Akh18} {\textsc{Akhmet, M., Fen, M. O.}, Non-autonomous equations with unpredictable solutions, \textit{Commun. Nonlinear Sci. Numer. Simulat.} \textbf{59} (2018), 657--670.}
					
				
		\bibitem{Akh19S} {\textsc{Akhmet, M., Fen, M. O., Tola, A.}, The sequential test for chaos, 2019, arXiv:1904.09127.}
		
		
		\bibitem{Casta} {\textsc{Casta\~{n}eda, L. B., Arunachalam, V., Dharmarajs, S.}, \textit{Introduction to Probability and Stochastic Processes with Applications}, Wiley, NJ, 2012.}	

		\bibitem{Wiggins} {\textsc{Wiggins, S.}, \textit{Global Bifurcation and Chaos: Analytical Methods},Springer-Verlag, New York, Berlin, 1988.}








	\end{thebibliography}
\end{document}